# Hidden variable recurrent fractal interpolation function with four function contractivity factors

Chol-Hui Yun

Faculty of Mathematics, **Kim Il Sung** University,

Pyongyang, Democratic People's Republic of Korea

Abstract: In this paper, we introduce a construction of hidden variable recurrent fractal interpolation functions (HVRFIF) with four function contractivity factors. In the fractal interpolation theory, it is very important to ensure flexibility and diversity of the construction of interpolation function. Recurrent iterated function system (RIFS) produce fractal sets with local self-similarity structure. Therefore the RIFS can describe the irregular and complicated objects in nature better than the iterated function system (IFS). Hidden variable fractal interpolation function (HVFIF) is neither self-similar nor self- affine one. The HVFIF is more complicated, diverse and irregular than the fractal interpolation function (FIF). The contractivity factor is important one that determins characteristics of FIFs. We present a constructions of one variable HVRFIFs and bivariable HVRFIFs using RIFS with four function contractivity factors.

**1. Introduction**

In 1986, Bsrnsley introduced a notation of fractal interpolation function(FIF) based on the theory of the iterated function system (IFS). FIFs have more advantages in modeling phenomena and fucntions with some self-similarity in nature than classical interpolation functions such as polynomials and spline. Therefore, FIFs have been studied in many articles and applied to a lot of areas such as function approxomation, signal process and computer graphics etc.

In fractal interpolation, one makes generally an IFS or a recurrent iterated function system (RIFS) for a given dataset and then constructs the Read-Bajraktarebic operator on suitable space of continuous functions using the IFS (RIFS). A fixed point of the operator is the FIF (RFIF) for the given dataset. Construction, derivative, integral, dimension, smoothness and stability of the FIFs have been widely studied. [1-16]

Since the vertical scaling factors, which the contraction transformations of IFS have, determine the characteristics of FIFs, they are very important. To obtain FIFs with high flexibility, construction of FIFs with function vertical scaling factors and their analytic properties have been studied in many papers. [6, 8, 12, 14, 16, 17]

RIFS is a generalization of IFS and produce local self-similar sets which aremore complicated than self-similar sets. FIF constructed by the RIFS is called a recurrent fractal interpolation function (RFIF). Constructions of RFIFs for a dataset in $R^2$ and $R^3$. [4, 5, 15] In [15], generalizing the construction of RFIFs with constant vertical scaling factors, RFIFs using the RIFS with function vertical scaling factors were constructed and the fractal dimension of graph of the constructed interpolation function was estimated and a construction of bivariate fractal interpolation functions using the RFIFs was proposed.

Barnsley *et al.* [2, 3] introduced a concept of hidden variable fractal interpolation function (HVFIF) which is more complicated, diverse and irregular than the FIF for the same set of interpolation data.

The idea of the construction of the HVFIF is to extend the given data set on $R^2$ into a data set on $R^3$, make a vector valued fractal interpolation function for the extended data set and then project the vector valued function onto $R^2$, which gives the HVFIF. It is usually non self-affine, because the HVFIF is the projection of a vector valued function. The HVFIFs have four free parameters: free variable, constrained free variable, free hidden variable and constrained free hidden variable (they are called contractivity factors). Using hidden variables, we can control more flexibly shapes and fractal dimensions of the graphs of FIFs.

In many papers, the contractivity factors are constants. [3, 6, 7] Therefore, the constructions lack the flexibility which is necessary to model complicated and irregular natural phenomena. To solve the problem, HVFIF with one function contractivity factor in [12] and HVFIF with four function contractivity factors were constructed.

In this paper, in order to ensure the flexibility and diversity of construction of HVFIF, we present construction of one variable and bivariable HVFIF with four function contractivity factors using RIFS. This paper is organized as follows:

In section 2, we construct one variable HVFIFs with four function contractivity factors for the extended data set. (Theorem 1, 2) In section 3, we present a construction of bivariable HVFIFs with four function contractivity factors (Theorem 3, 4).

## 2. One variable HVRFIF

### 2.1. Construction of RIFS

Let a data set $P_0$ in $R^2$ be given by

$$P_0 = \{(x_i, y_j) \in R^2; i = 0, 1, \cdots, n\}, \quad (-\infty < x_0 < x_1 < \cdots < x_n < +\infty). \quad (1).$$

To construct the HVFIF for the dataset $P_0$, we extend the dataset as follows:

$$P = \{(x_i, y_i, z_i,) = (x_i, \vec{y}_i) \in R^3; i = 0, 1, \cdots, n\}, \quad (-\infty < x_0 < x_1 < \cdots < x_n < \infty), \quad (2)$$

where $\vec{y}_i = (y_i, z_i)$ and $z_i$, $i = 0, 1, \cdots, n$ are parameters. Moreover, we denote $N_n = \{1, 2, \cdots, n\}$, $I_i = [x_{i-1}, x_i]$ and $I = [x_0, x_n]$, where $I_i$, $i \in \{1, \cdots, n\}$ is called a region. Let $l$ be integer with $2 \leq l \leq n$. We make subintervals $\tilde{I}_k$, $k = 1, \cdots, l$ of $I$ consisting of some regions and call $\tilde{I}_k$ a domain. Then, two end points of the domain $\tilde{I}_k$ ($k \in \{1, \cdots, l\}$) are contained in the set $\{x_0, x_1, \cdots, x_n\}$ and hence denoting start point and end point of $\tilde{I}_k$ by $s(k), e(k)$, respectively, we get the following mappings:

$$s: \{1, \cdots, l\} \to \{1, \cdots, n\}, \quad e: \{1, \cdots, l\} \to \{1, \cdots, n\}$$

and $\tilde{I}_k$ is denoted by $\tilde{I}_k = [x_{s(k)}, x_{e(k)}]$. We suppose that $e(k) - s(k) \geq 2$, $k = 1, \cdots, l$ which means that the interval $\tilde{I}_k$ contains at least 2 $I_i$ s.

For each $i \in N_n$, we take a $k (\in \{1, \cdots, l\})$ and denote it by $\gamma(i)$. Let mappings $L_{i,k} : [x_{s(k)}, x_{e(k)}] \to [x_{i-1}, x_i]$, $i \in N_n$ be contraction homeomorphisms that map end points of $\tilde{I}_k$ to end points of $I_i$, i.e. $L_{i,k}(\{x_{s(k)}, x_{e(k)}\}) = \{x_{i-1}, x_i\}$.

We define mappimgs $\vec{F}_{i,k} : \tilde{I}_k \times R^2 \to R^2$, $i = 1, \cdots, n$ as follows:

$$\vec{F}_{i,k}(x, \vec{y}) = \begin{pmatrix} s_i(L_{i,k}(x))y + s'_i(L_{i,k}(x))z + q_{i,k}(x) \\ \tilde{s}_i(L_{i,k}(x))y + \tilde{s}'_i(L_{i,k}(x))z + \tilde{q}_{i,k}(x) \end{pmatrix} = \begin{pmatrix} s_i(L_{i,k}(x)) & s'_i(L_{i,k}(x)) \\ \tilde{s}_i(L_{i,k}(x)) & \tilde{s}'_i(L_{i,k}(x)) \end{pmatrix} \begin{pmatrix} y \\ z \end{pmatrix} + \begin{pmatrix} q_{i,k}(x) \\ \tilde{q}_{i,k}(x) \end{pmatrix},$$

where $s_i, s'_i, \tilde{s}_i, \tilde{s}'_i : I_i \to R$ are Lipshitz functions on $I_i$ whose absolute value is less than 1 (which are called contractivity factors) and $q_{i,k}$, $\tilde{q}_{i,k} : \tilde{I}_k \to R$ are Lipshitz functions such that if $\alpha \in \{s(k), e(k)\}$, $L_{i,k}(x_\alpha) = x_a$ $a \in \{i-1, i\}$, then $\vec{F}_{i,k}(x_\alpha, \vec{y}_\alpha) = \vec{y}_a$. Then, $\vec{F}_{i,k}(x, \vec{z})$ is obviously Lipshitz function.

Let us denote $\vec{F}_{i,k}(x, \vec{y})$ by $\vec{F}_{i,k}(x, \vec{y}) = S_i \vec{y} + \vec{Q}_{i,k}(x)$, where

$$\vec{S}_i(x) = \begin{pmatrix} s_i(L_{i,k}(x)) & s_i'(L_{i,k}(x)) \\ \tilde{s}_i(L_{i,k}(x)) & \tilde{s}_i'(L_{i,k}(x)) \end{pmatrix}, \quad \vec{Q}_{i,k}(x) = \begin{pmatrix} q_{i,k}(x) \\ \tilde{q}_{i,k}(x) \end{pmatrix}, \quad \vec{y} = \begin{pmatrix} y \\ z \end{pmatrix}.$$

Example 1. An example of $q_{i,k}$, $\tilde{q}_{i,k}$ satisfied above is as follows:

$$q_{i,k}(x) = -s_i(L_{i,k}(x))g_{i,k}(x) - s_i'(L_{i,k}(x))g_{i,k}'(x) + h_i(L_{i,k}(x)),$$

$$\tilde{q}_{i,k}(x) = -\tilde{s}_i(L_{i,k}(x))g_{i,k}(x) - \tilde{s}_i'(L_{i,k}(x))g_{i,k}'(x) + \tilde{h}_i(L_{i,k}(x)),$$

$$g_{i,k}(x) = \frac{x - x_{s(k)}}{x_{e(k)} - x_{s(k)}} y_{e(k)} + \frac{x - x_{e(k)}}{x_{s(k)} - x_{e(k)}} y_{s(k)}, \quad g_{i,k}'(x) = \frac{x - x_{s(k)}}{x_{e(k)} - x_{s(k)}} z_{e(k)} + \frac{x - x_{e(k)}}{x_{s(k)} - x_{e(k)}} z_{s(k)},$$

$$h_i(x) = \frac{x - x_{i-1}}{x_i - x_{i-1}} y_i + \frac{x - x_i}{x_{i-1} - x_i} y_{i-1}, \quad \tilde{h}_i(x) = \frac{x - x_{i-1}}{x_i - x_{i-1}} z_i + \frac{x - x_i}{x_{i-1} - x_i} z_{i-1}.$$

Let $D \subset R^2$ be a sufficiently large bounded set containing $\vec{y}_i$, $i = 1, \cdots, n$.

We define transformations $\vec{W}_i : \tilde{I}_k \times D \to I_i \times R^2$, $i = 1, \cdots, n$ by

$$\vec{W}_i(x, \vec{y}) = (L_{i,k}(x), \vec{F}_{i,k}(x, \vec{y})), \quad i = 1, \cdots, n.$$

Then, as we know from the definitions of $L_{i,k}$ and $\vec{F}_{i,k}$, $\vec{W}_i$ maps the data points on end of the domain $\tilde{I}_k$ to the data points of $I_i$. For a function $f$, let us denote $\overline{f} = \max_x |f(x)|$. We denote $\overline{S} = \max\{\overline{s}_i + \overline{\tilde{s}}_i, \overline{s}_i' + \overline{\tilde{s}}_i'; i = 1, \cdots, n\}$. The following theorem gives a sufficient condition for $\vec{W}_i$ to be contraction transformation.

**Theorem 1.** If $\overline{S} < 1$, then there exists some distance $\rho_\theta$ equivalent to the Euclidean metric on $R^2$ such that $W_{i,k}$ ($i = 1, \cdots, n; k = \gamma(i)$) are contraction transformations with respect to the distance $\rho_\theta$.

Proof. We take $\theta$ as positive real number such that $\theta < \dfrac{1 - L_L}{L_s \alpha + L_Q}$, where

$$L_L = \max\{L_{L_i}; i = 1, \cdots, n\}, \quad L_S = \max\{L_{s_i}L_{L_{i,k}} + L_{\tilde{s}_i}L_{L_{i,k}}, L_{s_i'}L_{L_{i,k}} + L_{\tilde{s}_i'}L_{L_{i,k}}; i = 1, \cdots, n\},$$

$$L_Q = \max\{L_{q_{i,k}} + L_{\tilde{q}_{i,k}}; i = 1, \cdots, n\}, \quad \alpha = \sup_{\vec{y} \in D} \|\vec{y}\|_1 \text{ and } \|\cdot\|_1 \text{ is a norm on } R^2.$$

We define a distance $\rho_\theta$ on $R^3$ as follows:

$$\rho_\theta((x, \vec{y}), (x', \vec{y}')) = \|x - x'\|_1 + \theta \|\vec{y} - \vec{y}'\|_1, \quad (x, \vec{y}), (x', \vec{y}') \in R^3.$$

Then, it is clear that $\rho_\theta$ is equivalent to the Euclidean metric on $R^2$.

For $(x, \vec{y}), (x', \vec{y}') \in I \times D$, we have

$$\rho_\theta(\vec{W}_i(x, \vec{y}), \vec{W}_i(x', \vec{y}')) = \|\vec{L}_{i,k}(x) - \vec{L}_{i,k}(x')\|_1 + \theta \|\vec{F}_{i,k}(x, \vec{y}) - \vec{F}_{i,k}(x', \vec{y}')\|_1$$

$$= \|L_{i,k}(x) - L_{i,k}(x')\|_1 + \theta \|\vec{F}_{i,k}(x, \vec{y}) - \vec{F}_{i,k}(x', \vec{y}')\|_1$$

$$= \|L_{i,k}(x) - L_{i,k}(x')\|_1 + \theta \|\vec{S}_i(x)\vec{y} + \vec{Q}_{i,k}(x) - \vec{S}_i(x')\vec{y}' - \vec{Q}_{i,k}(x')\|_1$$

$$= \|L_{i,k}(x) - L_{i,k}(x')\|_1 + \theta \|(\vec{S}_i(x)\vec{y} - \vec{S}_i(x)\vec{y}') + (\vec{S}_i(x)\vec{y}' - \vec{S}_i(x')\vec{y}') + (\vec{Q}_{i,k}(x) - \vec{Q}_{i,k}(x'))\|_1$$

$$= \|L_{i,k}(x) - L_{i,k}(x')\|_1 + \theta \|\vec{S}_i(x)\|_1 \|\vec{y} - \vec{y}'\|_1 + \|\vec{S}_i(x) - \vec{S}_i(x')\|_1 \|\vec{y}'\|_1 + \|\vec{Q}_{i,k}(x) - \vec{Q}_{i,k}(x')\|_1$$

$$\leq L_L \|x - x'\|_1 + \overline{S}\theta \|\vec{y} - \vec{y}'\|_1 + L_S \alpha \|x - x'\|_1 + L_Q \theta \|x - x'\|_1$$

$$= (L_L + \theta(L_S \alpha + L_Q)) \|x - x'\| + \overline{S}\theta \|\vec{y} - \vec{y}'\|_1$$

$$\leq \max\{L_L + \theta(L_S\alpha + L_Q), \overline{S}\}(\|x - x'\|_1 + \theta\|\vec{y} - \vec{y}'\|_1)$$

$$= s\rho_\theta((x, \vec{y}), (x', \vec{y}')).$$

From the hypothesis of the theorem the condition on $\theta$, we have

$$s = \max\{L_L + \theta(L_S\alpha + L_Q), \overline{S}\} < 1.$$

Therefore, $\vec{W}_i, i = 1, \cdots, n$ are contraction transformations. □

We define a row-stochastic matrix $M = (p_{st})_{n \times n}$ by

$$p_{st} = \begin{cases} 1/a_s, & I_s \subseteq \tilde{I}_{\gamma(t)} \\ 0, & I_s \not\subset \tilde{I}_{\gamma(t)} \end{cases},$$

where for every $i(\in N_n)$, the number $a_i$ indicates the number of the domains $\tilde{I}_k$ containing the region $I_s$, which means that $p_{st}$ is positive if there is a transformation $W_i$ mapping $I_s$ to $I_t$.

Then, we have RIFS $\{R^3; M; W_i, i = 1, \cdots, n\}$ corresponding to the extended dataset **P**. We denote an attractor of the RIFS by A.

### 2.2 Construction of HVFIF

We make a continuous functuion that is a fixed point of the Read-Bajraktarebic operator defined by the RIFS, interpolates the extended dataset $P$ and whose graph is the attractor of the RIFS.

For the RIFS constructed above, we have the following theorem.

**Theorem 2. There is a continuous function $\vec{f}$ interpolating the extended data set $P$ such that the graph of $\vec{f}$ is the attractor A of RIFS constructed above.**

Proof. Let a set $\overline{C}(I)$ be as follows:

$\overline{C}(I) = \{\vec{h} : I \to R^2; \vec{h}$ interpolates the extended data set $P$ and is continuous$\}$

We can easily know that the set is a complete metric space with respect to the norm $\|\cdot\|_\infty$.

For $\vec{h}(\in \overline{C}(I))$, we define a mapping $T\vec{h}$ on $I$ by

$$(T\vec{h})(x) = \vec{F}_{i,k}(L_{i,k}^{-1}(x), \vec{h}(L_{i,k}^{-1}(x))), \quad x \in I_i.$$

Then, we get $T\vec{h} \in \overline{C}(I)$.

In fact, for any $i \in \{0, 1, \cdots, n\}$, there is $\alpha \in \{s(k), e(k)\}$ such that

$$(T\vec{h})(x_i) = \vec{F}_{i,k}(L_{i,k}^{-1}(x_i), \vec{h}(L_{i,k}^{-1}(x_i))) = \vec{F}_{i,k}(x_\alpha, \vec{y}_\alpha)) = \vec{y}_i, (L_{i,k}(x_\alpha) = x_i),$$

where $L_{0,k}(x) = L_{1,k}(x)$.

Hence, we can define an operator $T : \overline{C}(E) \to \overline{C}(E)$ on $\overline{C}(E)$. The operator is contraction one. In fact, we have

$$\|(T\vec{h})(x) - (T\vec{h}')(x)\|_1 = \|\vec{F}_{i,k}(L_{i,k}^{-1}(x), \vec{h}(L_{i,k}^{-1}(x))) - \vec{F}_{i,k}(L_{i,k}^{-1}(x), \vec{h}'(L_{i,k}^{-1}(x)))\|_1$$

$$= \|\vec{S}_i(x)\vec{h}(L_{i,k}^{-1}(x)) + \vec{Q}_{i,k}(L_{i,k}^{-1}(x)) - \vec{S}_i(x)\vec{h}'(L_{i,k}^{-1}(x)) - \vec{Q}_{i,k}(L_{i,k}^{-1}(x))\|_1$$

$$= \|\vec{S}_i(x)\|_1 \|\vec{h}'(L_{i,k}^{-1}(x)) - \vec{h}'(L_{i,k}^{-1}(x))\|_1$$

$$\leq \overline{S}\|h - h'\|_1.$$

Therefore, the operator $T$ has a unique fixed point $\vec{f}(\in \overline{C}(I))$ and

$$\vec{f}(x) = \vec{F}_{i,k}(L_{i,k}^{-1}(x), \vec{f}(L_{i,k}^{-1}(x))).$$

This gives that for the graph $Gr(\vec{f})$ of $\vec{f}$, $Gr(\vec{f}) = \bigcup_{j=1}^{n} \bigcup_{k \in I(j)} W_k(Gr(\vec{f}|_{E_k}))$,

where $I(j) = \{k \in \mathbb{N}; p_{kj} > 0\}$, $j = 1, \cdots, n$. This means that $Gr(\vec{f})$ is the attractor of the RIFS constructed in Theorem 1. Therefore, from the uniqueness of attractor, we have $A = Gr(\vec{f})$. □

Let us denote the vector valued function $\vec{f} : I \to \mathbb{R}^2$ in Theorem 2 by $\vec{f} = (f_1, f_2)$, where $f_1 : I \to \mathbb{R}$ interpolates the given dataset $P_0$, which is called a hidden variable recurrent fractal interpolation function (HVRFIF). Furthermore, a set $\{(x, f_1(x)) : x \in I\}$ is a projection of $A$ on $\mathbb{R}^2$.

As we know from the proof of Theorem 2, we have

$$\vec{f}(x) = \vec{F}_{i,k}(L_{i,k}^{-1}(x), \vec{f}(L_{i,k}^{-1}(x))), \quad x \in I_i.$$

i.e. $\vec{f}(x) = \vec{F}_{i,k}(L_{i,k}^{-1}(x), f_1(L_{i,k}^{-1}(x)), f_2(L_{i,k}^{-1}(x)))$, $x \in I_i$. Therefore, for all $x \in I$, the HVRFIF $f_1$ satisfies

$$f_1(x) = s_i(x) f_1(L_{i,k}^{-1}(x)) + s_i'(x) f_2(L_{i,k}^{-1}(x)) + q_{i,k}(L_{i,k}^{-1}(x))$$

and $f_2$ satisfies $f_2(x) = \tilde{s}_i(x) f_1(L_{i,k}^{-1}(x)) + \tilde{s}_i'(x) f_2(L_{i,k}^{-1}(x)) + \tilde{q}_{i,k}(L_{i,k}^{-1}(x))$.

Example 2. Figure 1 shows the graphs of one variable HVRFIFs constructed by RIFs with four function contractivity factors for a dataset

$P_0 = \{(0, 20), (0.25, 30), (0.5, 10), (0.75, 50), (1, 40)\}$.

Contractivity factors $\{s_1, s_2, s_3, s_4\}$, $\{\tilde{s}_1, \tilde{s}_2, \tilde{s}_3, \tilde{s}_4\}$, $\{s_1', s_2', s_3', s_{41}'\}$, $\{\tilde{s}_1', \tilde{s}_2', \tilde{s}_3', \tilde{s}_4'\}$ are as follows:

(1) $\{0.3, 0.85, 0.8, 0.5\}$, $\{0, 0, 0, 0\}$, $\{0.8, 0.6, 0.4, 0.5\}$, $\{0.19, 0.37, 0.48, 0.43\}$, (2) $\{0.3, 0.85, 0.8, 0.5\}$, $\{0.64, 0.14, 0.19, 0.49\}$, $\{0.8, 0.6, 0.4, 0.5\}$, $\{0.19, 0.37, 0.48, 0.43\}$, (3) $\{2.9x, 1.9x, x, x\}$, $\{0, 0, 0, 0\}$, $\{\sin(10x), \cos(300x), \sin(100x), \cos(3x)\}$, $\{0.99-|\sin(10x)|, 0.9-|\cos(300x)|, 0.95-|\sin(100x)|, 0.9-|\cos(3x)|\}$, (4) $\{2.9x, 1.9x, x, x\}$, $\{0.99-2.9x, 0.99-1.9x, 0.9-x, 0.99-x\}$, $\{\sin(10x), \cos(300x), \sin(100x), \cos(3x)\}$, $\{0.99-|\sin(10x)|, 0.9-|\cos(300x)|, 0.95-|\sin(100x)|, 0.9-|\cos(3x)|\}$.

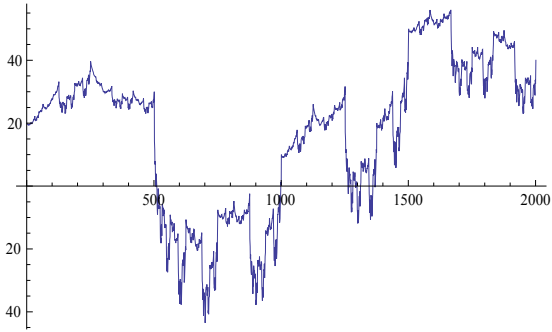
(1)

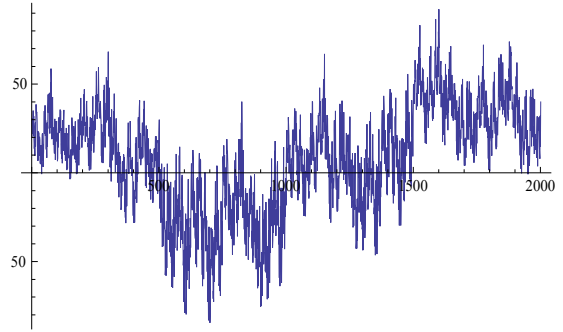
(2)

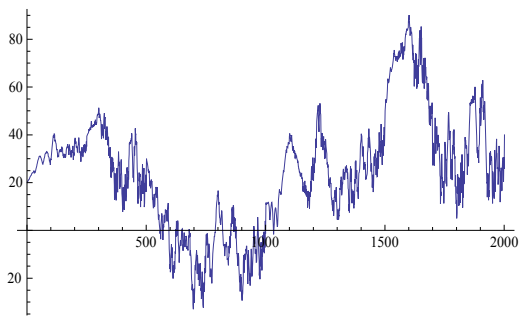
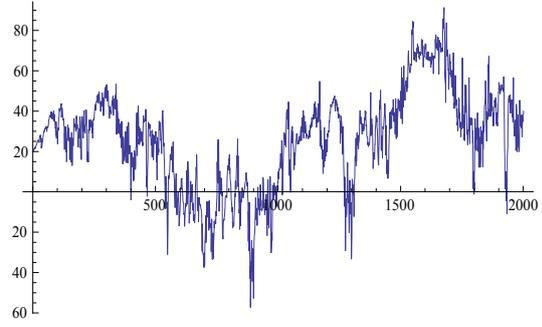

(3)           (4)

Fig.1. Graphs of HVRFIFs

## 3. Hidden variable bivariable recurrent FIF

**3.1.** Construction of RIFS

Let a dataset $P_0$ on rectangular grids be given as follows:
$$P_0 = \{(x_i, y_j, z_{ij}) \in R^3; i = 0, 1, \cdots, n, j = 0, 1, \cdots, m\}, \quad (x_0 < x_1 < \cdots < x_n, y_0 < y_1 < \cdots < y_m).$$
To make a HVRFIF for the dataset, we extend the dataset to the following one:
$$P = \{(x_i, y_j, z_{ij}, t_{ij}) = (\vec{x}_{ij}, \vec{z}_{ij}) \in R^4; i = 0,1,\cdots,n, j = 0,1,\cdots,m\}, \quad (x_0 < x_1 < \cdots < x_n, y_0 < y_1 < \cdots < y_m),$$
where $\vec{x}_{ij} = (x_i, y_j)$, $\vec{z}_{ij} = (z_{ij}, t_{ij})$ and $t_{ij}$, $i = 0, 1, \cdots, n, j = 0, 1, \cdots, m$ are parameters. We denote
$$N = n \cdot m, \ I_{x_i} = [x_{i-1}, x_i], \ N_{nm} = \{1,\cdots,n\} \times \{1,\cdots,m\}, \ I_{y_j} = [y_{j-1}, y_j], \ E = [x_0, x_n] \times [y_0, y_m], \ E_{ij} = I_{x_i} \times I_{y_j},$$
where $E_{ij}$ is called a region. Let $l$ be an integer with $2 \leq l \leq N$. Next, we take rectangulars $\widetilde{E}_k$, $k = 1,\cdots,l$ consisting of some regions from $E$. $\widetilde{E}_k$ is called a domain. Then we have $\widetilde{E}_k = \widetilde{I}_{x,k} \times \widetilde{I}_{y,k}$, where $\widetilde{I}_{x,k}$, $\widetilde{I}_{y,k}$ are closed intervals on x-axis and y-axis, respectively. Since the endpoints of $\widetilde{I}_{x,k}$ ($k \in \{1,\cdots,l\}$) are coincided with some endpoints of $I_{x_i}$, $i = 1,\cdots,n$, denoting start point and end point of $\widetilde{I}_{x,k}$ by $s_x(k), e_x(k)$, respectively, we can define the following mapping:
$$s_x : \{1,\cdots,l\} \to \{1,\cdots,n\}, \ e_x : \{1,\cdots,l\} \to \{1,\cdots,n\}.$$
Similarly, for $\widetilde{I}_{y,k}$, we define the mapping
$$s_y : \{1,\cdots,l\} \to \{1,\cdots,m\}, \ e_y : \{1,\cdots,l\} \to \{1,\cdots,m\}.$$
Then we have $\widetilde{I}_{x,k} = [x_{s_x(k)}, x_{e_x(k)}]$, $\widetilde{I}_{y,k} = [y_{s_y(k)}, y_{e_y(k)}]$, where we assume that $e_x(k) - s_x(k) \geq 2$, $e_y(k) - s_y(k) \geq 2$, $k = 1,\cdots,l$ which means that $\widetilde{I}_{x,k}$, $\widetilde{I}_{y,k}$ are intervals containing more than 2 $I_{x_i}$, $I_{x_j}$.

For $(i, j) \in N_{nm}$, we take $k (\in \{1,\cdots,l\})$ and denote it by $\gamma(i, j)$.

We define mappings
$$L_{x_i,k} : [x_{s_x(k)}, x_{e_x(k)}] \to [x_{i-1}, x_i], \ L_{y_j,k} : [y_{s_y(k)}, y_{e_y(k)}] \to [y_{j-1}, y_j], \ (i, j) \in N_{nm}$$
as contraction homeomorphisms that map end points of $\widetilde{I}_{x,k}$, $\widetilde{I}_{y,k}$ to end points of $I_{x_i}$, $I_{y_j}$, i.e.
$$L_{x_i,k}(\{x_{s_x(k)}, x_{e_x(k)}\}) = \{x_{i-1}, x_i\}, \ L_{y_j,k}(\{y_{s_y(k)}, y_{e_y(k)}\}) = \{y_{j-1}, y_j\}.$$
Next, we define transformations $\vec{L}_{ij} : \widetilde{E}_k \to E_{ij}$ by $\vec{L}_{ij}(\vec{x}) = (L_{x_i,k}(x), L_{y_j,k}(y))$. Then these $\vec{L}_{ij}$ map vertexes of $\widetilde{E}_k$ to one of $E_{ij}$, i.e. for $\alpha \in \{s_x(k), e_x(k)\}$, $\beta \in \{s_y(k), e_y(k)\}$,
$$\vec{L}_{ij}(\vec{x}_{\alpha\beta}) = \vec{x}_{ab} \ (a \in \{i-1, i\}, \ b \in \{j-1, j\}).$$

We define mappings $\vec{F}_{ij} : \widetilde{E}_k \times R^2 \to R^2$, $i = 1,\cdots,n, \ j = 1,\cdots,m$ by
$$\vec{F}_{ij,k}(\vec{x}, \vec{z}) = \begin{pmatrix} s_{ij}(\vec{L}_{ij}(\vec{x}))z + s'_{ij}(\vec{L}_{ij}(\vec{x}))t + q_{ij}(\vec{x}) \\ \widetilde{s}_{ij}(\vec{L}_{ij}(\vec{x}))z + \widetilde{s}'_{ij}(\vec{L}_{ij}(\vec{x}))t + \widetilde{q}_{ij}(\vec{x}) \end{pmatrix} = \begin{pmatrix} s_{ij}(\vec{L}_{ij}(\vec{x})) & s'_{ij}(\vec{L}_{ij}(\vec{x})) \\ \widetilde{s}_{ij}(\vec{L}_{ij}(\vec{x})) & \widetilde{s}'_{ij}(\vec{L}_{ij}(\vec{x})) \end{pmatrix} \begin{pmatrix} z \\ t \end{pmatrix} + \begin{pmatrix} q_{ij}(\vec{x}) \\ \widetilde{q}_{ij}(\vec{x}) \end{pmatrix},$$
where $s_{ij}, s'_{ij}, \widetilde{s}_{ij}, \widetilde{s}'_{ij} : E_{ij} \to R$ are arbitrary Lipschitz functions whose absolute values are less than 1 and $q_{ij}, \widetilde{q}_{ij} : \widetilde{E}_k \to R$ are defined as mappings satisfying the following condition: for

$\alpha \in \{s_x(k), e_x(k)\}$, $\beta \in \{s_y(k), e_y(k)\}$, $L_{x_i,k}(x_\alpha) = x_a$, $L_{y_j,k}(y_\beta) = y_b$ ($a \in \{i-1, i\}$, $b \in \{j-1, j\}$),

$$\vec{F}_{ij,k}(\vec{x}_{\alpha\beta}, \vec{z}_{\alpha\beta}) = \vec{z}_{ab}.$$

Then, $\vec{F}_{ij,k}(\vec{x}, \vec{z})$ are Lipschitz mappings.

We denote $\vec{F}_{ij,k}(\vec{x}, \vec{z})$ by $\vec{F}_{ij,k}(\vec{x}, \vec{z}) = \vec{S}_{ij}\vec{z} + \vec{Q}_{ij}(\vec{x})$, where

$$\vec{S}_{ij}(\vec{x}) = \begin{pmatrix} s_{ij}(L_{ij}(\vec{x})) & s'_{ij}(L_{ij}(\vec{x})) \\ \tilde{s}_{ij}(L_{ij}(\vec{x})) & \tilde{s}'_{ij}(L_{ij}(\vec{x})) \end{pmatrix}, \quad \vec{Q}_{ij}(\vec{x}) = \begin{pmatrix} q_{ij}(\vec{x}) \\ \tilde{q}_{ij}(\vec{x}) \end{pmatrix}, \quad \vec{z} = \begin{pmatrix} z \\ t \end{pmatrix}.$$

In the future, we denote simply $\vec{F}_{ij,k}$, $L_{x_i,k}$, $L_{y_j,k}$ by $\vec{F}_{ij}$, $L_{x_i,k}$, $L_{y_j,k}$. Let $D \subset R^2$ be a sufficiently large bounded set containing $\vec{z}_{ij}$, $i = 1, \cdots, n$, $j = 1, \cdots, m$.

Now, we define transformation $\vec{W}_{ij} : \tilde{E}_k \times D \to E_{ij} \times R^2$, $i = 1, \cdots, n$, $j = 1, \cdots, m$ by

$$\vec{W}_{ij}(\vec{x}, \vec{z}) = (\vec{L}_{ij}(\vec{x}), \vec{F}_{ij}(\vec{x}, \vec{z})), \quad i = 1, \cdots, n, \, j = 1, \cdots, m$$

We denote $\overline{S} = \max\{\bar{s}_{ij} + \bar{\tilde{s}}_{ij}, \bar{s}'_{ij} + \bar{\tilde{s}}'_{ij}; i = 1, \cdots, n, j = 1, \cdots, m\}$. The following theorem gives a sufficient condition for $\vec{W}_{ij}$ to be contraction one.

**Theorem 3.** If $\overline{S} < 1$, then there exists some distance $\rho_\theta$ equivalent to the Euclidean metric such that $\vec{W}_{ij}$, $i = 1, \cdots, n$, $j = 1, \cdots, m$ are contraction transformations with respect to the distance $\rho_\theta$.

Proof. We take $\theta$ as positive number such that

$$\theta < \frac{1 - c_L}{L_S \alpha + L_Q}, \tag{3}$$

where $c_{\vec{L}_{ij}} = \max\{c_{L_{x_i}}, c_{Ly_j}\}$, $c_L = \max\{c_{\vec{L}_{ij}}; i = 1, \cdots, n, j = 1, \cdots, m\}$,

$L_S = \max\{L_{s_{ij}} c_{\vec{L}_{ij}} + L_{\tilde{s}_{ij}} c_{\vec{L}_{ij}}, L_{s'_{ij}} c_{\vec{L}_{ij}} + L_{\tilde{s}'_{ij}} c_{\vec{L}_{ij}}; i = 1, \cdots, n, j = 1, \cdots, m\}$,

$L_Q = \max\{L_{q_{ij}} + L_{\tilde{q}_{ij}}; i = 1, \cdots, n, j = 1, \cdots, m\}$, $\alpha = \sup_{\vec{z} \in D} \|\vec{z}\|_1$,

and $\|\cdot\|_1$ is a norm on $R^2$.

Now, we define a distance $\rho_\theta$ on $R^4$ by

$$\rho_\theta((\vec{x}, \vec{z}), (\vec{x}', \vec{z}')) = \|\vec{x} - \vec{x}'\|_1 + \theta \|\vec{z} - \vec{z}'\|_1, \quad (\vec{x}, \vec{z}), (\vec{x}', \vec{z}') \in R^4.$$

Then, $\rho_\theta$ is obviously equivalent to the Euclidean metric on $R^4$ and for $(\vec{x}, \vec{z}), (\vec{x}', \vec{z}') \in E \times D$, we have

$\rho_\theta(\vec{W}_{ij}(\vec{x}, \vec{z}), \vec{W}_{ij}(\vec{x}', \vec{z}')) =$

$= \|\vec{L}_{ij}(\vec{x}) - \vec{L}_{ij}(\vec{x}')\|_1 + \theta \|\vec{F}_{ij}(\vec{x}, \vec{z}) - \vec{F}_{ij}(\vec{x}', \vec{z}')\|_1$

$= \|\vec{L}_{ij}(\vec{x}) - \vec{L}_{ij}(\vec{x}')\|_1 + \theta \|\vec{S}_{ij}(\vec{x})\vec{z} + \vec{Q}_{ij}(\vec{x}) - \vec{S}_{ij}(\vec{x}')\vec{z}' - \vec{Q}_{ij}(\vec{x}')\|_1$

$= \|\vec{L}_{ij}(\vec{x}) - \vec{L}_{ij}(\vec{x}')\|_1 + \theta \|(\vec{S}_{ij}(\vec{x})\vec{z} - \vec{S}_{ij}(\vec{x})\vec{z}') + (\vec{S}_{ij}(\vec{x})\vec{z}' - \vec{S}_{ij}(\vec{x}')\vec{z}') + (\vec{Q}_{ij}(\vec{x}) - \vec{Q}_{ij}(\vec{x}'))\|_1$

$= \|\vec{L}_{ij}(\vec{x}) - \vec{L}_{ij}(\vec{x}')\|_1 + \theta \|\vec{S}_{ij}(\vec{x})\|_1 \|\vec{z} - \vec{z}'\|_1 + \|\vec{S}_{ij}(\vec{x}) - \vec{S}_{ij}(\vec{x}')\|_1 \|\vec{z}'\|_1 + \|\vec{Q}_{ij}(\vec{x}) - \vec{Q}_{ij}(\vec{x}')\|_1$

$\leq c_L \|\vec{x} - \vec{x}'\|_1 + \overline{S}\theta \|\vec{z} - \vec{z}'\|_1 + L_S \alpha \|\vec{x} - \vec{x}'\|_1 + L_Q \theta \|\vec{x} - \vec{x}'\|_1$

$= (c_L + \theta(L_S \alpha + L_Q)) \|\vec{x} - \vec{x}'\| + \overline{S}\theta \|\vec{z} - \vec{z}'\|_1$

$\leq \max\{c_L + \theta(L_S \alpha + L_Q), \overline{S}\}(\|\vec{x} - \vec{x}'\|_1 + \theta \|\vec{z} - \vec{z}'\|_1)$

$= s\rho_\theta((\vec{x}, \vec{z}), (\vec{x}', \vec{z}'))$.

From the hypothesis of the theorem and condition on $\theta$, we have $s = \max\{c_L + \theta(L_S \alpha + L_Q), \overline{S}\} < 1$.

This means that $\vec{W}_{ij}, i = 1, \cdots, n, j = 1, \cdots, m$ are contraction transformations. □

Remark: In the definition of $\rho_\theta$, even in the case where $\|\cdot\|_1$ is changed into $\|\cdot\|_\infty$, we have the similar result.

We define a row-stochastic matrix $M = (p_{st})_{N \times N}$ by

$$p_{st} = \begin{cases} 1/a_s, & E_{\tau^{-1}(s)} \subseteq \widetilde{E}_{\gamma(\tau^{-1}(t))} \\ 0, & E_{\tau^{-1}(s)} \not\subset \widetilde{E}_{\gamma(\tau^{-1}(t))} \end{cases},$$

where $\tau : N_{nm} \to \{1, \cdots, N\}$ is an one to one mapping defined by $\tau(i,j) = i + (j-1)n$ and the number $a_s$ indicates the number of the domains $\widetilde{E}_k$, $k = 1, \cdots, l$ containing the region $E_{\tau^{-1}(s)}$, which means that $p_{st}$ is positive if there is a transformation $W_{ij}$ mapping $E_s$ to $E_t$.

Then, RIFS $\{R^4; M; \vec{W}_{ij}, i = 1, \cdots, n, j = 1, \cdots, m\}$ is a RIFS corresponding to the extended dataset P. An attractor of the RIFS is called a recurrent fractal set and denoted by $A$.

**2) Construction of bivariable HVRFIFs**

We present a sufficient condition for a fixed point of the Read-Bajraktarebic operator defined by the RIFS to interpolate the extended dataset $P$ and have a graph which is the attractor of the RIFS.

In the mapping $\vec{F}_{ij,k}(\vec{x}, \vec{z}) = \vec{S}_{ij}\vec{z} + \vec{Q}_{ij}(\vec{x})$ defined above, we define $\vec{Q}_{ij}$ as one satisfying condition that for some continuous function $\vec{g}$ interpolating the dataset $P$, i.e. $\vec{g}(\vec{x}_{ij}) = \vec{z}_{ij}$, $i = 0, 1, \cdots, n$, $j = 0, 1, \cdots, m$,

$$\vec{F}_{ij}(x_\alpha, y, \vec{g}(x_\alpha, y)) = \vec{g}(\vec{L}_{ij}(x_\alpha, y)), \ \alpha \in \{s_x(k), e_x(k)\}, \tag{3}$$

$$\vec{F}_{ij}(x, y_\beta, \vec{g}(x, y_\beta)) = \vec{g}(\vec{L}_{ij}(x, y_\beta)), \ \beta \in \{s_y(k), e_y(k)\}. \tag{4}$$

One example is as follows:

$$\vec{Q}_{ij}(\vec{x}) = -\vec{S}_{ij}(\vec{L}_{ij}(\vec{x}))\vec{l}_{ij}(\vec{x}) + \vec{r}_{ij}(\vec{L}_{ij}(\vec{x})),$$

where

$$\vec{r}_{ij}(\vec{x}) = \vec{g}(\vec{x}), \ \vec{x} \in \partial E_{ij},$$

$$\vec{l}_{ij}(\vec{x}) = \vec{g}(\vec{x}), \ \vec{x} \in \partial \widetilde{E}_{\gamma(i,j)},$$

Then, we have $\vec{F}_{ij}(\vec{x}, \vec{z}) = \vec{S}_{ij}(\vec{L}_{ij}(\vec{x}))(\vec{z} - \vec{l}_{ij}(\vec{x})) + \vec{r}_{ij}(\vec{L}_{ij}(\vec{x}))$. Next example is one of $\vec{l}_{ij}(\vec{x})$ and $\vec{r}_{ij}(\vec{x})$ satisfying these conditions

Example 3. A mapping $\vec{g}(x, y)$ is an interpolation function of the given dataset. $\vec{l}_{ij}(\vec{x})$ and $\vec{r}_{ij}(\vec{x})$ coincide with $\vec{g}|_{\partial \widetilde{E}_k}$, $\vec{g}|_{\partial E_{ij}}$ on the domain $\partial \widetilde{E}_k$ the region $\partial E_{ij}$, respectively. For instance, we define

$$\vec{l}_{ij}(\vec{x}) := \vec{g}(x_i, y)(1 - u_l) + \vec{g}(x_{i+2}, y)u_l + \vec{g}(x, y_j)(1 - v_l) + \vec{g}(x, y_{j+2})v_l -$$
$$- (1 - u_l)(1 - v_l)\vec{z}_{ij} - (1 - u_l)v_l \vec{z}_{ij+2} - u_l(1 - v_l)\vec{z}_{i+2j} - u_l v_l \vec{z}_{i+2j+2}$$

$$u_l := \frac{x - x_i}{x_{i+2} - x_i}, \quad v_l := \frac{y - y_j}{y_{j+2} - y_j}, \quad (x, y) \in [x_i, x_{i+2}] \times [y_j, y_{j+2}]$$

$$\vec{r}_{ij}(x, y) := \vec{g}(x_i, y)(1 - u_r) + \vec{g}(x_{i+1}, y)u_r + \vec{g}(x, y_j)(1 - v_r) + \vec{g}(x, y_{j+1})v_r -$$
$$- (1 - u_r)(1 - v_r)\vec{z}_{ij} - (1 - u_r)v_r \vec{z}_{ij+1} - u_r(1 - v_r)\vec{z}_{i+1j} - u_r v_r \vec{z}_{i+1j+1}$$

$$u_r := \frac{x - x_i}{x_{i+1} - x_i}, \quad v_r := \frac{y - y_j}{y_{j+1} - y_j}, \quad (x, y) \in [x_i, x_{i+1}] \times [y_j, y_{j+1}]$$

Then, we have

$\vec{l}_{ij}(x_\alpha, y) = \vec{g}(x_\alpha, y)$, $\vec{l}_{ij}(x, y_\beta) = \vec{g}(x, y_\beta)$, $\alpha \in \{s_x(k), e_x(k)\}$, $\beta \in \{s_y(k), e_y(k)\}$, $\vec{r}_{ij}(x_a, y) = \vec{g}(x_a, y)$, $\vec{r}_{ij}(x, y_b) = g(x, y_b)$, $a \in \{i-1, i\}$, $b \in \{j-1, j\}$).

Let the attractor of the RIFS satisfying conditions (3) and (4) by B, we have the following theorem.

**Theorem 4. There is a continuous function $\vec{f}$ which interpolates the data set P and whose graph is the attractor B.**

Proof. We define a set $\overline{C}(E)$ by

$$\overline{C}(E) = \{\vec{h} : E \to R^2 \in C; \vec{h} \text{ interpolates } P \text{ and satisfies conditions (3), (4)}\}.$$

We can easily prove that this is a complete metric space with respect to the norm $\|\cdot\|_\infty$. Since $\vec{g} \in \overline{C}(E)$, which is contained in the definition of $\vec{F}_{ij}$, we have $\overline{C}(E) \neq \varnothing$. Moreover, the functions coinciding with $\vec{g}$ on $\partial E_{ij}$ are contained in $\overline{C}(E)$.

For $\vec{h}(\in \overline{C}(E))$, we define a mapping $T\vec{h}$ on $E$ by
$$(T\vec{h})(x, y) = \vec{F}_{ij}(\vec{L}_{ij}^{-1}(x, y), \vec{h}(\vec{L}_{ij}^{-1}(x, y))), \quad (x, y) \in E_{ij}.$$

Then, we have $T\vec{h} \in \overline{C}(E)$.

In fact, for $\alpha \in \{s_x(k), e_x(k)\}$,
$$(T\vec{h})(\vec{L}_{ij}(x_\alpha, y)) = \vec{F}_{ij}(x_\alpha, y, \vec{h}(x_\alpha, y)) = \vec{h}(\vec{L}_{ij}(x_\alpha, y)),$$
$$(T\vec{h})(\vec{L}_{ij}(x, y_\beta)) = \vec{F}_{ij}(x, y_\beta, \vec{h}(x, y_\beta)) = \vec{h}(\vec{L}_{ij}(x, y_\beta)).$$

Therefore, $T\vec{h} = \vec{h}$ on $\{(x_i, y): y \in [y_0, y_m]\}$, $\{(x, y_j): x \in [x_0, x_n]\}$, $i = 0, 1, \cdots, n$, $j = 0, 1, \cdots, m$ and
$$\vec{F}_{ij}(x_\alpha, y, (T\vec{h})(x_\alpha, y)) = \vec{F}_{ij}(x_\alpha, y, \vec{h}(x_\alpha, y)) = \vec{h}(L_{ij}(x_\alpha, y)) = (T\vec{h})(\vec{L}_{ij}(x_\alpha, y))$$

i.e. $T\vec{h}$ satisfies (3). Similarly, we can prove that $T\vec{h}$ satisfies (4).

Therefore, we can define an operator $T: \overline{C}(E) \to \overline{C}(E)$ on $\overline{C}(E)$. It is obvious that the operator is contraction one. In fact, we have

$$\|(T\vec{h})(x, y) - (T\vec{h}')(x, y)\|_1 = \|\vec{F}_{ij}(\vec{L}_{ij}^{-1}(x, y), \vec{h}(\vec{L}_{ij}^{-1}(x, y))) - \vec{F}_{ij}(\vec{L}_{ij}^{-1}(x, y), \vec{h}'(\vec{L}_{ij}^{-1}(x, y)))\|_1$$
$$= \|\vec{S}_{ij}(x, y)\vec{h}(\vec{L}_{ij}^{-1}(x, y)) + \vec{Q}_{ij}(\vec{L}_{ij}^{-1}(x, y)) - \vec{S}_{ij}(x, y)\vec{h}'(\vec{L}_{ij}^{-1}(x, y)) - \vec{Q}_{ij}(\vec{L}_{ij}^{-1}(x, y))\|_1$$
$$= \|\vec{S}_{ij}(x, y)\|_1 \|\vec{h}'(\vec{L}_{ij}^{-1}(x, y)) - h'(\vec{L}_{ij}^{-1}(x, y))\|_1$$
$$\leq \overline{S} \|h - h'\|_1.$$

Hence, $T$ has a unique fixed point $\vec{f}(\in \overline{C}(E))$,
$$\vec{f}(x, y) = \vec{F}_{ij}(\vec{L}_{ij}^{-1}(x, y), \vec{f}(\vec{L}_{ij}^{-1}(x, y))).$$

Furthermore, for the graph $Gr(\vec{f})$ of $\vec{f}$, we get $Gr(\vec{f}) = \bigcup_{j=1}^{N} \bigcup_{k \in I(j)} W_{\tau^{-1}(k)}(Gr(\vec{f}|_{E_{\tau^{-1}(k)}}))$.

This means that $Gr(\vec{f})$ is the attractor of RIFS constructed above. From the uniqueness of the attractor of RIFS, we have $A = Gr(\vec{f})$. □

The vector valued function $\vec{f}: E \to R^2$ in Theorem 4 is denoted by $\vec{f} = (f_1(\vec{x}), f_2(\vec{x}))$, where $f_1: E \to R$ interpolates the given dataset $P_0$ and is called hidden variable bivariable recurrent fractal interpolation function (HVBRFIF) for the dataset $P_0$. Moreover, $\{(\vec{x}, f_1(\vec{x})): \vec{x} \in E\}$ is a projection of B on $R^3$. Since the projection is not always self-affine, the HVBRFIF is not generally self-affine FIF. $f_2(\vec{x})$ interpolates the set $\{(x_i, y_j, t_{ij}) = (\vec{x}_{ij}, t_{ij}) \in R^3; i = 0, 1, \cdots, n, j = 0, 1, \cdots, m\}$.

From the proof of Theorem 4, we get
$$\vec{f}(x, y) = \vec{F}_{ij}(\vec{L}_{ij}^{-1}(x, y), \vec{f}(\vec{L}_{ij}^{-1}(x, y))), \quad (x, y) \in E_{ij},$$

i.e. $\vec{f}(x, y) = \vec{F}_{ij}(\vec{L}_{ij}^{-1}(x, y), f_1(\vec{L}_{ij}^{-1}(x, y)), f_2(\vec{L}_{ij}^{-1}(x, y)))$, $(x, y) \in E_{ij}$.

Therefore, for all $(x, y) \in E$, HVBRFIF $f_1$ satisfies
$$f_1(x, y) = s_{ij}(x, y) f_1(\vec{L}_{ij}^{-1}(x, y)) + s'_{ij}(x, y) f_2(\vec{L}_{ij}^{-1}(x, y)) + q_{ij}(\vec{L}_{ij}^{-1}(x, y))$$

and $f_2$ satisfies $f_2(x, y) = \tilde{s}_{ij}(x, y) f_1(\vec{L}_{ij}^{-1}(x, y)) + \tilde{s}'_{ij}(x, y) f_2(\vec{L}_{ij}^{-1}(x, y)) + \tilde{q}_{ij}(\vec{L}_{ij}^{-1}(x, y))$.

Example 2. Figure 2 shows the graphs of HVBRFIFs constructed by RIFSs with some contractivity factors and dataset $P_0$ given in the following table:

| y \ x | 0 | 0.25 | 0.5 | 0.75 | 1 |
|---|---|---|---|---|---|
| 0 | 46 | 32 | 65 | 73 | 39 |
| 0.25 | 32 | 23 | 84 | 33 | 29 |
| 0.5 | 76 | 88 | 58 | 73 | 88 |
| 0.75 | 62 | 79 | 33 | 86 | 43 |
| 1 | 49 | 23 | 39 | 76 | 32 |

Contractivity factors $\{s_{ij}\}$, $\{\tilde{s}_{ij}\}$, $\{s'_{ij}\}$ and $\{\tilde{s}'_{ij}\}$ are as follows:

(1)

|   | 1 | 2 | 3 | 4 |
|---|---|---|---|---|
| 1 | 0.9 | 0.6 | -0.9 | 0.7 |
| 2 | -0.94 | 0.95 | 0.3 | 0.85 |
| 3 | 0.5 | -0.99 | 0.86 | 0.79 |
| 4 | 0.87 | 0.92 | 0.75 | -0.87 |

|   | 1 | 2 | 3 | 4 |
|---|---|---|---|---|
| 1 | -0.4 | 0.9 | -0.65 | 0.66 |
| 2 | 0.9 | 0.36 | 0.27 | 0.25 |
| 3 | 0.91 | -0.89 | 0.9 | 0.85 |
| 4 | 0.53 | 0.96 | -0.49 | 0.39 |

|   | 1 | 2 | 3 | 4 |
|---|---|---|---|---|
| 1 | 0 | 0 | 0 | 0 |
| 2 | 0 | 0 | 0 | 0 |
| 3 | 0 | 0 | 0 | 0 |
| 4 | 0 | 0 | 0 | 0 |

|   | 1 | 2 | 3 | 4 |
|---|---|---|---|---|
| 1 | 0.53 | 0.03 | 0.27 | 0.28 |
| 2 | 0.09 | 0.55 | 0.64 | 0.72 |
| 3 | 0.05 | 0.01 | 0.02 | 0.11 |
| 4 | 0.41 | 0.03 | 0.48 | 0.56 |

,

(2)

|   | 1 | 2 | 3 | 4 |
|---|---|---|---|---|
| 1 | 0.9 | 0.6 | -0.9 | 0.7 |
| 2 | -0.94 | 0.95 | 0.3 | 0.85 |
| 3 | 0.5 | -0.99 | 0.86 | 0.79 |
| 4 | 0.87 | 0.92 | 0.75 | -0.87 |

|   | 1 | 2 | 3 | 4 |
|---|---|---|---|---|
| 1 | -0.4 | 0.9 | -0.65 | 0.66 |
| 2 | 0.9 | 0.36 | 0.27 | 0.25 |
| 3 | 0.91 | -0.89 | 0.9 | 0.85 |
| 4 | 0.53 | 0.96 | -0.49 | 0.39 |

|   | 1 | 2 | 3 | 4 |
|---|---|---|---|---|
| 1 | -0.47 | -0.07 | -0.08 | 0.14 |
| 2 | 0.15 | -0.04 | -0.69 | 0.14 |
| 3 | 0.46 | -0.69 | 0.04 | 0.07 |
| 4 | 0.07 | -0.13 | 0.18 | -0.02 |

|   | 1 | 2 | 3 | 4 |
|---|---|---|---|---|
| 1 | 0.53 | 0.03 | 0.27 | 0.28 |
| 2 | 0.09 | 0.55 | 0.64 | 0.72 |
| 3 | 0.05 | 0.01 | 0.02 | 0.11 |
| 4 | 0.41 | 0.03 | 0.48 | 0.56 |

(3)

|   | 1 | 2 | 3 |
|---|---|---|---|
| 1 | 0.45(cos(x)+sin(y)) | 0.9sin(10x^2+10y^2) | 0.9cos(10x^2+10y^2) |
| 2 | 0.99cos(10x^3+10y^3) | 0.45(cos(x)-sin(y)) | 0.9sin(x^2+y^9) |
| 3 | 0.9cos(50x+50y) | 0.99cos(40x^3+40y^3) | 0.9sin(10x+10y) |
| 4 | 0.9sin(50x+50y) | 0.99cos(150x^2+15y^2) | 0.45(cos(20x)-sin(20y)) |

|   | 1 | 2 | 3 |
|---|---|---|---|
| 1 | 0 | 0 | 0 |
| 2 | 0 | 0 | 0 |
| 3 | 0 | 0 | 0 |
| 4 | 0 | 0 | 0 |

|   | 1 | 2 | 3 |
|---|---|---|---|
| 1 | 0.45(cos(x)+sin(y)) | 0.9sin(10x^2+10y^2) | 0.9cos(30x^2+30y^2) |
| 2 | 0.93cos(10x^4+10y^3) | 0.45(cos(x)-sin(y)) | 0.95sin(x^2+y^9) |
| 3 | 0.9cos(20x+20y) | 0.85cos(30x^3+40y^3) | 0.87sin(20x+30y) |
| 4 | 0.98sin(40x+40y) | 0.93cos(150x^2+15y^2) | 0.4(cos(30x)-sin(20y)) |

|   | 1 | 2 | 3 |
|---|---|---|---|
| 1 | 0.93-|0.45(cos(x)+sin(y))| | 0.93-|0.9sin(10x^2+10y^2)| | 0.92-|0.9cos(30x^2+30y^2)| |
| 2 | 0.99-|0.93cos(10x^4+10y^3)| | 0.91-|0.45(cos(x)-sin(y))| | 0.91-|0.95sin(x^2+y^9)| |
| 3 | 0.96-|0.9cos(20x+20y)| | 0.99-|0.85cos(30x^3+40y^3)| | 0.92-|0.87sin(20x+30y)| |
| 4 | 0.94-|0.98sin(40x+40y)| | 0.99-|0.93cos(150x^2+15y^2)| | 0.97-|0.4(cos(30x)-sin(20y))| |

(4)

|   | 1 | 2 | 3 |
|---|---|---|---|
| 1 | 0.45(cos(x)+sin(y)) | 0.9sin(10x^2+10y^2) | 0.9cos(10x^2+10y^2) |
| 2 | 0.99cos(10x^3+10y^3) | 0.45(cos(x)-sin(y)) | 0.9sin(x^2+y^9) |
| 3 | 0.9cos(50x+50y) | 0.99cos(40x^3+40y^3) | 0.9sin(10x+10y) |
| 4 | 0.9sin(50x+50y) | 0.99cos(150x^2+15y^2) | 0.45(cos(20x)-sin(20y)) |

|   | 1 | 2 | 3 |
|---|---|---|---|
| 1 | 0.82-|0.45(cos(x)+sin(y))| | 0.84-|0.9sin(10x^2+10y^2)| | 0.82-|0.9cos(10x^2+10y^2|) |
| 2 | 0.79-|0.99cos(10x^3+10y^3)| | 0.91-|0.45(cos(x)-sin(y))| | 0.99-|0.9sin(x^2+y^9)| |
| 3 | 0.96-|0.9cos(50x+50y)| | 0.69-|0.99cos(40x^3+40y^3)| | 0.9-|0.9sin(10x+10y)| |
| 4 | 0.94-|0.9sin(50x+50y)| | 0.79-|0.99cos(150x^2+15y^2)| | 0.93-|0.45(cos(20x)-sin(20y))| |

|   | 1 | 2 | 3 |
|---|---|---|---|
| 1 | 0.45(cos(x)+sin(y)) | 0.9sin(10x^2+10y^2) | 0.9cos(30x^2+30y^2) |
| 2 | 0.93cos(10x^4+10y^3) | 0.45(cos(x)-sin(y)) | 0.95sin(x^2+y^9) |
| 3 | 0.9cos(20x+20y) | 0.85cos(30x^3+40y^3) | 0.87sin(20x+30y) |
| 4 | 0.98sin(40x+40y) | 0.93cos(150x^2+15y^2) | 0.4(cos(30x)-sin(20y)) |

|   | 1 | 2 | 3 |
|---|---|---|---|
| 1 | 0.93-|0.45(cos(x)+sin(y))| | 0.93-|0.9sin(10x^2+10y^2)| | 0.92-|0.9cos(30x^2+30y^2)| |
| 2 | 0.99-|0.93cos(10x^4+10y^3)| | 0.91-|0.45(cos(x)-sin(y))| | 0.91-|0.95sin(x^2+y^9)| |
| 3 | 0.96-|0.9cos(20x+20y)| | 0.99-|0.85cos(30x^3+40y^3)| | 0.92-|0.87sin(20x+30y)| |
| 4 | 0.94-|0.98sin(40x+40y)| | 0.99-|0.93cos(150x^2+15y^2)| | 0.97-|0.4(cos(30x)-sin(20y))| |

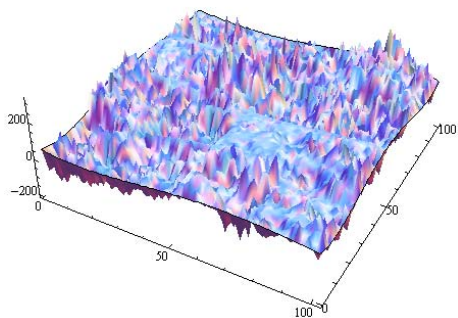
(1)

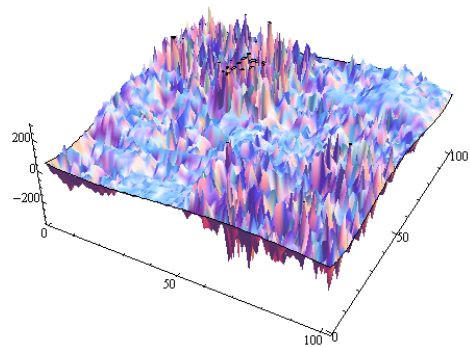
(2)

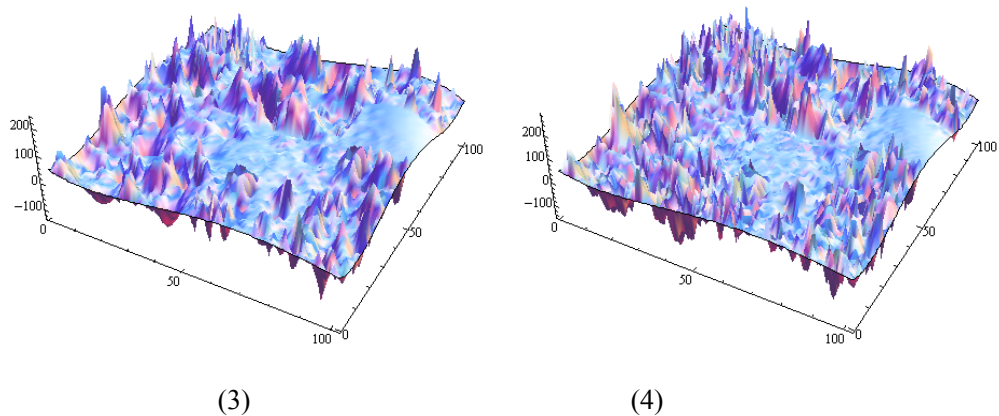

(3) (4)

Fig.2. Graphs of HVBRFIFs